\documentclass[12pt]{article}
\usepackage{graphicx,pst-node,pst-plot,pstricks,fancybox}
\usepackage{latexsym,amssymb,amsmath}
\usepackage{times}
\usepackage{locenglish,math}
\usepackage{graphicx,pst-plot,psfig,pstricks}
\begin{document}
\begin{center}
\textbf{\large{\textsf{Random walks on randomly oriented lattices}}\footnote{
\textit{1991 Mathematics Subject Classification:}
60J10, 60K20\\
\textit{Key words and phrases:}
Markov chain, random environment, recurrence criteria, random graphs, oriented graphs.
}}
\vskip1cm
\parbox[t]{14cm}{
Massimo {\sc  Campanino}$^a$ and Dimitri {\sc Petritis}$^b$\\
\vskip5mm   
{\scriptsize  
\baselineskip=5mm
a. Dipartimento di Matematica, Università degli Studi di Bologna,\\
piazza di Porta San Donato 5, I-40126 Bologna, Italy, campanin@dm.unibo.it
\vskip5mm
b. Institut de Recherche
Math\'ematique, Universit\'e de Rennes I and CNRS UMR 6625\\
Campus de Beaulieu, F-35042 Rennes Cedex, France, Dimitri.Petritis@univ-rennes1.fr}\\
\vskip1cm
{\small
\centerline{19 November 2001}
\vskip1cm
\baselineskip=5mm
{\bf Abstract:}

Simple random walks on various types of partially horizontally oriented regular lattices are considered.
The horizontal orientations of the lattices can be of various types (deterministic or random)
 and depending on the
nature of the orientation the asymptotic behaviour of the random walk is shown to be recurrent
or transient. In particular, for randomly horizontally oriented lattices the random walk is almost
surely transient. 
}}
\end{center}

\section{Introduction}
\subsection{Motivations}
Random walks are mathematical objects with important applications in many
scientific disciplines and in particular in physics.
Although the bulk of this paper is devoted to the probabilistic problems
arising for a particular class of random walks,
some indications on the physical
interest of the objects we introduce will be given briefly in this subsection.

Beyond the original impetus for the study of random walks 
given by the seminal work
of Einstein on diffusion ---  
an informal but fascinating account of which can be found
in chapter 5 of \cite{Pai}, ---
there was a revival of the physical interest 
for the subject in the early '80 because it
allowed a rigorous and powerful representation of Green's
 functions in Euclidean (scalar) quantum field theory
and statistical mechanics (see \cite{FerFroSok} for a review).
This representation serves also as a rigorous basis for the numerical
simulation of quantities of physical relevance in
those two theories that remain otherwise inaccessible 
by the analytic computation.
With respect to this latter aspect, the denumerable graph on which the
random walk evolves is a discretised approximation of the continuum
space(-time) manifold.
The main drawback of the random walks on lattices 
is that they don't allow the study of
quantum field theories  more complicated 
than the scalar ones, like the gauge field theories
or the fermionic field theories, or of quantum statistical physics
 because these 
theories are intrinsically non-commutative
even in their Euclidean version. The discretised differential calculus for
these theories becomes the study of differential forms on 
graded algebras, necessitating
thus the  introduction of oriented 
(directed) lattices as discretised
versions of the space(-time) continuum \cite{DimMue}. 

Although random walks on oriented lattices 
are the relevant objects to study in
the context of discretised gauge theories, 
their rigorous probabilistic study is
still lacking. To the best of the authors knowledge, 
the only prior probabilistic
work on the topic is a paragraph containing a side-result in the PhD thesis
\cite{Guillotin}.

From a purely mathematical point of view, random walks on directed lattices
present also very interesting features. For instance, simple
random walks on undirected regular lattices (like $\BbZ^d$) are
thoroughly studied and a vast literature establishes precise criteria
for their transience or null recurrence properties. 
Not to mention but one result,
the recurrence of the  simple random walk on an undirected lattice  
is related to the convergence norm $r(\bP)$ of the transition matrix 
$\bP$ and the
latter is determined by the geometric properties of the graph
through its transition isoperimetric number 
(see corollary 5.6 of \cite{MohWoe}).
We shall see that this is not any longer the case for directed  graphs
since we shall exhibit two different  
regular deterministic directed graphs having the
same isoperimetric number the one  
being null recurrent and the other transient.
This remark constitutes our main  
motivation for studying random walks on randomly oriented
lattices. Since the choice of different deterministic orientation of the graph
leaves enough room to have dramatically different asymptotic behaviour for the
simple random walk defined on them, a natural question is how a random choice
of orientations would affect the result.

\subsection{Notations and definitions}
An \textit{oriented} (or equivalently \textit{directed}) graph 
$\BbG=(\BbV, \BbA)$ is the pair of a denumerable
set 
$\BbV$  of vertices and a set 
$\BbA\subset \BbV\times\BbV$
of oriented edges. We exclude the presence of loops 
(\textit{i.e.}\ edges $a=(\bv,\bv)$ with $\bv\in \BbV$).
Multiple edges are also excluded by definition. 
The corresponding graph is then termed \textit{simple}.

 \textit{Range} and a \textit{source} functions, denoted respectively $r$ and $s$,
are defined as mappings $r,s:\BbA\rightarrow \BbV$, defined by 
$\BbA\ni a=(\bu,\bv)\mapsto r(a)=\bv\in\BbV$ and
$\BbA\ni a=(\bu,\bv)\mapsto s(a)=\bu\in\BbV$.
We can therefore define, for each vertex $\bv\in\BbV$, its 
\textit{inwards degree}
$d^+_\bv=\card\{a\in\BbA: r(a)=\bv\}$ and its \textit{outwards degree}
$d^-_\bv=\card\{a\in\BbA: s(a)=\bv\}$.
All the graphs we consider are \textit{finitely transitive} in the sense
for any two distinct vertices $\bu,\bv\in\BbV$, there is a \textit{finite}
sequence $(\bw_0,\ldots,\bw_k)$ of vertices $\bw_i\in\BbV$, for $i=0,\ldots,k$,
$k\in\BbN^*$, with $\bw_0=\bu$ and $\bw_k=\bv$, 
such that $(\bw_i,\bw_{i+1})\in\BbA, \forall i=0,\ldots, k-1$. This proprerty
implies in particular the  \textit{no sink condition}:
$d^-_\bv\geq 1$ for all $\bv\in\BbV$.
Notice that undirected graphs can be considered as directed ones
verifying the condition that whenever an edge $(\bu, \bv)\in \BbA$ then
the reverse edge $(\bv,\bu)\in\BbA$. Therefore, when we speak about
directed graphs in the sequel, we mean general graphs where some edges
can be non-directed. However, we always consider
graphs that  are genuinely oriented in the sense that
there exist vertices $\bu$ and $\bv$ with
$(\bu,\bv)\in\BbA$ but $(\bv,\bu)\not\in\BbA$.

\begin{defi}{[Simple random walk]} Let
$(\BbV,\BbA)$ be an oriented graph. 
A \textit{simple random walk} on $(\BbV,\BbA)$
is a $\BbV$-valued Markov chain 
$(\bM_n)_{n\in\BbN}$ with transition probability matrix $\bP$
having as  matrix elements
\[P(\bu,\bv)=\BbP(\bM_{n+1}=\bv|\bM_n=\bu)=
\left\{\begin{array}{ll}
\frac{1}{d^-_\bu} & \textrm{if }\ \ (\bu,\bv)\in\BbA\\
0 & \textrm{otherwise.}
\end{array}\right.\]
\end{defi}

\Rk When the underlying graph is genuinely oriented, the Markov chain 
$(\bM_n)_{n\in\BbN}$ \textit{cannot be reversible}. Therefore, all the powerful techniques
based on the analogy with electrical circuits (see \cite{DoySne,Soa} for instance)
do not apply.
Notice moreover that the Markov operator of this chain \textit{is not expressible
in terms of the Laplace-Beltrami operator}; in the oriented case the Markov operator
is expressible in terms of the Dirac operator of the lattice (roughly, 
the square root
of the Laplace-Beltrami operator).

All the graphs that we shall consider in this paper are
two-dimensional lattices, i.e.\ $\BbV=\BbZ^2$ 
and $\BbA$ is a subset of the set of nearest neighbours in $\BbZ^2$.
We often write 
$\BbV=\BbV_1\times\BbV_2$, with $\BbV_1$ and $\BbV_2$ isomorphic to $\BbZ$
when we wish to specify horizontal and vertical directions.
In the latter notation, when necessary, 
we can distinguish between abscissas and ordinates
of vertices $\bv\in\BbV$ by writing $\bv=(v_1,v_2)$.

Let $\bom{\epsilon}=(\epsilon_y)_{y\in\BbV_2}$ 
be a $\{-1,1\}$-valued sequence of variables assigned to each
ordinate. The sequence $\bom{\epsilon}$ can be 
deterministic or random as it will
be specified later.

\begin{defi}{[$\bom{\epsilon}$-horizontally oriented lattice]} 
Let $\BbV=\BbV_1\times\BbV_2=\BbZ^2$, with $\BbV_1$ and $\BbV_2$ isomorphic to $\BbZ$
and $\bom{\epsilon}=(\epsilon_y)_{y\in\BbV_2}$ be a sequence of 
$\{-1,1\}$-valued variables assigned to each
ordinate. We call $\bom{\epsilon}$-\textit{horizontally oriented lattice} 
$\BbG=\BbG(\BbV,\bom{\epsilon})$,
the directed graph with vertex set $\BbV=\BbZ^2$ and edge set $\BbA$
defined by the condition
$(\bu,\bv)\in\BbA$ if, and only if, $\bu$ and $\bv$ are distinct vertices
satisfying one of the following conditions:
\begin{enumerate}
\item
either $v_1=u_1$ and $v_2=u_2\pm1$, 
\item
or $v_2=u_2$ and $v_1=u_1+\epsilon_{u_2}$.
\end{enumerate}
\end{defi}
\Rk Notice that the $\bom{\epsilon}$-horizontally oriented lattice
is regular; this means that the vertex degrees (both inwards and outwards)
are constant $d^-_\bv=d^+_\bv=d=3$, $\forall \bv\in\BbV$. 
The vertical directions of the graph are both-ways; the horizontal
directions are one-way, the sign of $\epsilon_y$ determining whether
the horizontal line at level $y$ is left- or right-going.

\begin{exam}{[Alternate lattice $\BbL$]}
In that case, $\bom{\epsilon}$ is the deterministic sequence $\epsilon_y=(-1)^y$ for
$y \in\BbV_2$.
The figure \ref{fig-alternate} depicts a part of this graph.
\end{exam}
\begin{figure}[h]
\centerline{%
\hbox{%
\psset{unit=4.5mm}
\pspicture(-0.5,-0.5)(7,7)
\Cartesian(4.5mm,4.5mm)
\psline[linewidth=0.1pt]{<-<}(-0.5,0)(6.5,0) 
\psline[linewidth=0.1pt]{>->}(-0.5,1)(6.5,1)
\psline[linewidth=0.1pt]{<-<}(-0.5,2)(6.5,2)
\psline[linewidth=0.1pt]{>->}(-0.5,3)(6.5,3)
\psline[linewidth=0.1pt]{<-<}(-0.5,4)(6.5,4)
\psline[linewidth=0.1pt]{>->}(-0.5,5)(6.5,5)   
\psline[linewidth=0.1pt]{<-<}(-0.5,6)(6.5,6)
\psline[linewidth=0.1pt]{-}(0,-0.5)(0,6.5)
\psline[linewidth=0.1pt]{-}(1,-0.5)(1,6.5)
\psline[linewidth=0.1pt]{-}(2,-0.5)(2,6.5)
\psline[linewidth=0.1pt]{-}(3,-0.5)(3,6.5)
\psline[linewidth=0.1pt]{-}(4,-0.5)(4,6.5)
\psline[linewidth=0.1pt]{-}(5,-0.5)(5,6.5)
\psline[linewidth=0.1pt]{-}(6,-0.5)(6,6.5)
\rput[bl](2.5,3.2){\tiny{$(0,0)$}}
\endpspicture}}
\caption{\label{fig-alternate}
\textsf{The alternately directed lattice $\BbL$ corresponding to the choice
 $\epsilon_{y}=(-1)^{y}$.}}
\end{figure}

\begin{exam}{[The half-plane one-way lattice $\BbH$]}
Here $\bom{\epsilon}$ is the deterministic sequence
\[\epsilon_y=\left\{\begin{array}{ll}
1 & \textrm{if } \ \ y\geq 0\\
-1 & \textrm{if } \ \ y<0.
\end{array}
\right.\]
The figure \ref{fig-halfplane} depicts a part of this graph.
\end{exam}

\begin{figure}[h]
\centerline{%
\hbox{%
\psset{unit=4.5mm}
\pspicture(-0.5,-0.5)(7,7)
\Cartesian(4.5mm,4.5mm)
\psline[linewidth=0.1pt]{<-<}(-0.5,0)(6.5,0)
\psline[linewidth=0.1pt]{<-<}(-0.5,1)(6.5,1)
\psline[linewidth=0.1pt]{<-<}(-0.5,2)(6.5,2)
\psline[linewidth=0.1pt]{>->}(-0.5,3)(6.5,3)
\psline[linewidth=0.1pt]{>->}(-0.5,4)(6.5,4)
\psline[linewidth=0.1pt]{>->}(-0.5,5)(6.5,5)
\psline[linewidth=0.1pt]{>->}(-0.5,6)(6.5,6)
\psline[linewidth=0.1pt]{-}(0,-0.5)(0,6.5)
\psline[linewidth=0.1pt]{-}(1,-0.5)(1,6.5)  
\psline[linewidth=0.1pt]{-}(2,-0.5)(2,6.5)  
\psline[linewidth=0.1pt]{-}(3,-0.5)(3,6.5)  
\psline[linewidth=0.1pt]{-}(4,-0.5)(4,6.5)  
\psline[linewidth=0.1pt]{-}(5,-0.5)(5,6.5)  
\psline[linewidth=0.1pt]{-}(6,-0.5)(6,6.5)  
\rput[bl](2.5,3.2){\tiny{$(0,0)$}}
\endpspicture}}
\caption{\label{fig-halfplane}
\textsf{The half-plane one-way  lattice $\BbH$ with
$\epsilon_y=-1$, if $y<0$ and
$\epsilon_y=1$, if $y\geq0$.}}
\end{figure}

\begin{exam}{[The lattice with random horizontal 
orientations $\BbO_{\bom{\epsilon}}$]}
Here $\bom{\epsilon}=(\epsilon_y)_{y\in\BbV_2}$ 
is a sequence of Rademacher, \textit{i.e.}\
 $\{-1,1\}$-valued 
symmetric Bernoulli random variables,
that are independent for different values of $y$.
The figure \ref{fig-randomly} depicts part of a realisation of this graph.
The random sequence  $\bom{\epsilon}$ is also termed the 
\textit{environment of random horizontal
directions}.
\end{exam}

\begin{figure}[h]
\centerline{%
\hbox{%
\psset{unit=4.5mm}
\pspicture(-0.5,-0.5)(7,7)
\Cartesian(4.5mm,4.5mm)
\psline[linewidth=0.1pt]{<-<}(-0.5,0)(6.5,0)
\psline[linewidth=0.1pt]{<-<}(-0.5,1)(6.5,1)
\psline[linewidth=0.1pt]{>->}(-0.5,2)(6.5,2)
\psline[linewidth=0.1pt]{<-<}(-0.5,3)(6.5,3)
\psline[linewidth=0.1pt]{>->}(-0.5,4)(6.5,4)
\psline[linewidth=0.1pt]{>->}(-0.5,5)(6.5,5)
\psline[linewidth=0.1pt]{<-<}(-0.5,6)(6.5,6)
\psline[linewidth=0.1pt]{-}(0,-0.5)(0,6.5)
\psline[linewidth=0.1pt]{-}(1,-0.5)(1,6.5)
\psline[linewidth=0.1pt]{-}(2,-0.5)(2,6.5)
\psline[linewidth=0.1pt]{-}(3,-0.5)(3,6.5)
\psline[linewidth=0.1pt]{-}(4,-0.5)(4,6.5)
\psline[linewidth=0.1pt]{-}(5,-0.5)(5,6.5)
\psline[linewidth=0.1pt]{-}(6,-0.5)(6,6.5)
\rput[bl](2.5,3.2){\tiny{$(0,0)$}}
\endpspicture}}
\caption{\label{fig-randomly}
\textsf{The randomly horizontally 
directed  lattice $\BbO_{\bom{\epsilon}}$ with
$(\epsilon_y)_{y\in\BbZ}$ an independent and identically distributed
sequence of Rademacher random variables.}}
\end{figure}

\subsection{Results}
The graphs defined previously  are topologically non-trivial in the sense that
\[\lim\limits_{N\rightarrow\infty}\frac{1}{N}\sum_{y=-N}^N \epsilon_y=0.\]
For the lattices $\BbL$ and $\BbH$ this is shown by a simple calculation
and for the lattice $\BbO_{\bom{\epsilon}}$ this is an almost sure statement
stemming from the independence of the sequence $\bom{\epsilon}$.
The above condition guarantees the existence of 
infinitely many non trivial allowed loops
having the origin as one of their vertices.

We   are now ready to state our results.

\begin{theo}
\label{th-L}
The simple random walk on the alternate lattice $\BbL$ is recurrent.
\end{theo}

\Rk This result can be easily generalised to any lattice
with  periodically alternating  horizontal directions (for every
finite period).

\begin{theo}
\label{th-H}
The simple random walk on the half-plane one-way lattice $\BbH$ is transient.
\end{theo}

\Rk The result concerning transience in theorem \ref{th-H} is robust.  
In particular, perturbing the orientation of any finite set of
horizontal lines either by reversing the orientation of these
lines or by transforming them into two-ways does not change the
transient behaviour of the simple random walk.
Therefore, the half-plane one-way lattice is so deeply in the
transience region that the asymptotic behaviour of the
simple random walk cannot be changed by 
simply modifying the transition probabilities along
a lower dimensional manifold as  was 
the case in \cite{MenPet-wedge} where the bulk
behaviour is on the critical point and it can 
be changed by lower-dimensional perturbations.

\begin{theo}
\label{th-O}
For almost all realisations of the environment $\bom{\epsilon}$, 
the simple random
walk on the randomly horizontally oriented lattice 
$\BbO_{\bom{\epsilon}}$ is transient and its speed is $0$.
\end{theo}


\section{Technical preliminaries}
\subsection{Embedding}
We suppose that there is an abstract probability space
$(\Omega,\cA,\BbP)$ on which are defined  
all the random variables we shall use.
In particular, the Markov chain $(\bM_n)_{n\in\BbN}$
 is defined on this space
and we denote $(\tilde{\cF_n})_{n\in\BbN}$ the natural filtration
of the process  $(\bM_n)_{n\in\BbN}$.
We assume a deterministic anchoring of the chain
at the origin, \textit{i.e.}\ $\bM_0=(0,0)\in\BbV$; obviously
 $\tilde{\cF}_0\subseteq\ldots\tilde{\cF}_n\subseteq
\tilde{\cF}_{n+1}\subseteq\ldots  \tilde{\cF_\infty}\subseteq\cA$. 
Let $\be_1$ and $\be_2$ denote the unit vectors of $\BbZ^2$.
We denote $(\tilde{\psi}_n)_{n\in\BbN}$
the sequence of the vertical projections of the increments of the
Markov chain, more precisely,
\[
\tilde{\psi}_{n+1}= (\be_2, \bM_{n+1}-\bM_n),\]
where $(\cdot,\cdot)$ denotes the Euclidean scalar product on $\BbR^2$ and
$\be_1,\be_2$ the canonical basis of $\BbR^2$.
The random variables 
$(\tilde{\psi}_n)_{n\in\BbN,y\in\BbV_2}$
form a 
sequence of independent $\{-1,0,1\}$-valued
random variables, symmetrically distributed
according to the law
\begin{eqnarray*}
\BbP(\tilde{\psi}_i=-1)&=& 
\BbP(\tilde{\psi}_i=1)\ =\ \frac{1}{d^-}\\
\BbP(\tilde{\psi}_i= 0)&=& 1-\frac{2}{d^-}= p,
\end{eqnarray*}
where $d^-=3$ is the outwards degree of any vertex (recall that the
graphs we consider are regular) and $q=1-p=2/3$ represents
the probability that the walk performs a vertical move.
\begin{lemm}
On the event
$\{\bM_n=\bu\}$, the Markov chain $(\bM_n)$ verifies
\[\bM_{n+1}=\bu+\epsilon_{u_2}\be_1\id_{\{\tilde{\psi}_{n+1}=0\}}+
\tilde{\psi}_{n+1}\be_2\id_{\{\tilde{\psi}_{n+1}\not=0\}}.\]
\end{lemm}
\Proof Obvious. \eproof

Introduce now the infinite sequence of random times 
$(\tilde{\sigma}_k, \tilde{\tau}_k)_{k\in\BbN^*}$ by
\begin{eqnarray*}
\tilde{\sigma}_1&=&1\\
\tilde{\tau}_1	&=&\inf\{n\geq \tilde{\sigma}_1:
	\tilde{\psi}_n\ne 0\}\geq \tilde{\sigma}_1\\
\tilde{\sigma}_2&=&\inf\{n\geq \tilde{\tau}_1:
	\tilde{\psi}_n =  0\}> \tilde{\tau}_1\\
  &\vdots & 					     \\
\tilde{\tau}_k	&=&\inf\{n\geq \tilde{\sigma}_{k}:
	\tilde{\psi}_n\ne 0\}>\tilde{\sigma}_{k}\\
\tilde{\sigma}_{k+1}&=&\inf\{n\geq \tilde{\tau}_k:
	\tilde{\psi}_n =  0\}>\tilde{\tau}_k\\
&\vdots & 					     \\
\end{eqnarray*}

The idea of the embedding is to decompose the two-dimensional
random walk $(\bM_n)$ into a skeleton simple one-dimensional random walk 
--- the vertical one $(Y_n)$ ---, a sequence of waiting times
$(\xi_n)$, and
an embedded
one-dimensional random walk with unbounded jumps --- the horizontal 
one $(X_n)$. In order to achieve the decomposition of the random walk, 
regroup the instants $n$ into blocks separated by the symbol $|$ as follows:
\[\tilde{\sigma}_1\ldots \tilde{\tau}_1-1\ |\ 
\tilde{\tau}_1\ldots\tilde{\sigma}_2-1\ |\
\tilde{\sigma}_2 \ldots \tilde{\tau}_2-1\ |\ 
\tilde{\tau}_2\ldots \tilde{\sigma}_3-1\ |\
\tilde{\sigma}_3\ldots \tilde{\tau}_3-1\ |\ \ldots\]

For the instants 
in the blocks starting with a $\tilde{\sigma}$ --- notice that
 the leftmost
one may be empty but all the other contain at least one instant ---
the walk performs horizontal moves, for blocks
starting with a $\tilde{\tau}$ the walk performs vertical moves.
More precisely, 
define a sequence of random sets of integers for $k\in\BbN^*$ by
\begin{eqnarray*}
I_k &=& [\tilde{\sigma}_k, \tilde{\tau}_k-1]\cap \BbN\\
J_k &=& [\tilde{\tau}_k, \tilde{\sigma}_{k+1}-1]\cap \BbN,
\end{eqnarray*}
with $I_1$ being the empty set when $\tilde{\tau}_1=\tilde{\sigma}_1$.
The random walk $(\bM_n)$ performs horizontal moves when $n$ is in a
$I_k$ for some $k\in\BbN^*$ and vertical moves when $n$ is in a $J_k$
for some $k\in\BbN^*$. 

Shrink now the $I$ sets and replace them by a waiting time. More
precisely, if $I_1\not=\emptyset$, define $\alpha=1$ and $\tilde{\xi}_1=
\tilde{\tau}_1-\tilde{\sigma}_1=|I_1|$, where $|\cdot|$ denotes cardinality,
else define $\alpha=0$ and $\tilde{\xi}_0$ need not be defined.
Then recursively, for $n\in\BbN^*$, define
\begin{eqnarray}
\label{random-intervals}
\tilde{\xi}_{\sum_{i=1}^n|J_i|+\alpha+(n-1)}&=&
\tilde{\tau}_{n+1}-\tilde{\sigma}_{n},\nonumber\\
\tilde{\xi}_{\sum_{i=1}^n|J_i|+\alpha+(n-1)+k}&=& 
0,\ \ \forall k\in\{1,\ldots,|J_{n+1}|\},\\
\psi_{\sum_{i=1}^n|J_i|+k}&=& 
\tilde{\psi}_{\sum_{i=1}^n(|I_i|+|J_i|)+k}\in \{-1,1\},
\ \ \forall k\in\{1,\ldots,|J_{n+1}|\}.\nonumber
\end{eqnarray}
Notice that $\tilde{\xi}_l$ cannot be non-zero for two consecutive
indices since $|J_{n+1}|=\tilde{\sigma}_{n+2}-\tilde{\tau}_{n+2}\geq 1$
for every $n=1,2,\ldots$.

\begin{lemm}
Given a realisation of the sequence $(\tilde{\psi}_n)_{n\in\BbN^*}$
the sequences $(\tilde{\xi}_k)_{k\in\BbN^*}$ and 
$(\psi_k)_{k\in\BbN^*}$ are uniquely determined and conversely.
\end{lemm}
\Proof
The previous construction proves the direct way of the lemma.
To prove the converse, given the sequences
$(\tilde{\xi}_k)_{k\in\BbN^*}$ and 
$(\psi_k)_{k\in\BbN^*}$, the sequence $(\tilde{\psi}_n)_{n\in\BbN^*}$
is obtained 
by inflating the time $\tilde{\xi}_k$ spent in waiting 
to reconstruct the intervals $(I_k, J_k)_{k\in\BbN^*}$.
Then invert the relations (\ref{random-intervals})
by assigning the value $\tilde{\psi}_i=0$ whenever
$i$ belongs to an interval of $I$-type.
\eproof

The figure (\ref{fig-reform}) depicts an example of random walk on
the lattice and the table (\ref{tab-reform}) establishes the bijection
between the various sequences.

\begin{figure}[h]
\centerline{%
\hbox{%
\psset{unit=4.5mm}
\pspicture(-0.5,-0.5)(7,7)
\Cartesian(4.5mm,4.5mm)
\psline[linewidth=0.1pt]{<-<}(-0.5,0)(6.5,0)
\psline[linewidth=0.1pt]{<-<}(-0.5,1)(6.5,1)
\psline[linewidth=0.1pt]{>->}(-0.5,2)(6.5,2)
\psline[linewidth=0.1pt]{<-<}(-0.5,3)(6.5,3)
\psline[linewidth=0.1pt]{>->}(-0.5,4)(6.5,4)
\psline[linewidth=0.1pt]{>->}(-0.5,5)(6.5,5)
\psline[linewidth=0.1pt]{<-<}(-0.5,6)(6.5,6)
\psline[linewidth=0.1pt]{-}(0,-0.5)(0,6.5)
\psline[linewidth=0.1pt]{-}(1,-0.5)(1,6.5)
\psline[linewidth=0.1pt]{-}(2,-0.5)(2,6.5)
\psline[linewidth=0.1pt]{-}(3,-0.5)(3,6.5) 
\psline[linewidth=0.1pt]{-}(4,-0.5)(4,6.5)
\psline[linewidth=0.1pt]{-}(5,-0.5)(5,6.5)
\psline[linewidth=0.1pt]{-}(6,-0.5)(6,6.5)
\psline[linewidth=1pt]{-}(3,3)(3,4)
\psline[linewidth=1pt]{-}(3,4)(5,4)
\psline[linewidth=1pt]{-}(5,4)(5,1)
\psline[linewidth=1pt]{-}(5,1)(1,1)
\psline[linewidth=1pt]{-}(1,1)(1,2)
\psline[linewidth=1pt]{-}(1,2)(4,2)
\psline[linewidth=1pt]{->}(4,2)(4,3)
\rput[br](3.5,2.5){\tiny{$(0,0)$}}
\endpspicture}}
\caption{\label{fig-reform} 
\textsf{A realisation of the Markov chain $(\bM_n)_{n=0,\ldots, 15}$.}}
\end{figure}

\begin{table}[h]
\caption{\label{tab-reform} 
\textsf{The reformulation of the random walk according to vertical and horizontal moves.}}

\begin{center}
\begin{tabular}{||c|c|c|c||}
\hline  
$n$ 	& $\bM$		& $\psi$ 	& $\tilde{\xi}$ \\
\hline
$0$	& $(0,0)$ 	& 		& $0$	\\
$1$	& $(0,1)$	& $1$		& $2$  	\\
$2$     & $(1,1)$       &               &       \\
$3$     & $(2,1)$       &               &       \\
$4$     & $(2,0)$       & $-1$          & $0$   \\
$5$     & $(2,-1)$      & $-1$          & $0$   \\    
$6$     & $(2,-2)$      & $-1$          & $4$   \\ 
$7$     & $(1,-2)$      &               &       \\
$8$     & $(0,-2)$      &               &       \\
$9$     & $(-1,-2)$     &               &       \\
$10$    & $(-2,-2)$     &               &       \\ 
$11$    & $(-2,-1)$     & $1$		& $3$	\\
$12$    & $(-1,-1)$     &               &       \\ 
$13$    & $(0,-1)$      &               &       \\
$14$    & $(1,-1)$      &               &       \\
$15$    & $(1,0)$       & $1$           &       \\
\hline
\end{tabular}
\end{center}
\end{table}

\begin{lemm}
The sequence $(\psi_n)_{n\in\BbN^*}$ is an independent identically distributed
sequence of symmetric Bernoulli $\{-1,1\}$-valued random variables.
\end{lemm}
\Proof
The  independence follows from the independence of the 
$(\tilde{\psi}_n)_{n\in\BbN^*}$  sequence. For every $n\in\BbN^*$, the law
of  $\psi_n$ is the conditional law of a $\tilde{\psi}_m$ with respect
to the  event $\{\tilde{\psi}_m\not=0\}$.
\eproof

\begin{lemm}
The sequence $(\tilde{\xi}_n)_{n\in\BbN^*}$
is an independent, identically distributed sequence of $\BbN$-valued
geometric random
variables of parameters $p$ and $q=1-p$ with
\[\BbP(\tilde{\xi}_1=\ell)=pq^\ell, \ \ \ell=0,1,2,\ldots.\]
\end{lemm}
\Proof
The  independence follows from the independence of the 
$(\tilde{\psi}_n)_{n\in\BbN^*}$  sequence. For every $m\in\BbN^*$, 
the variable $\tilde{\xi}_m$ is nothing else than the waiting time
on the state $0$ for the sequence $(\tilde{\psi}_n)_{n\in\BbN^*}$.
\eproof

\subsection{Basic definitions}

\begin{defi}{}
Let $(\psi_n)_{n\in\BbN^*}$ be a sequence of independent, identically 
distributed, $\{-1,1\}$-valued symmetric Bernoulli variables and 
\[Y_n=\sum_{k=1}^n\psi_k, \ \ n=1,2,\ldots\]
with $Y_0$, the simple $\BbV_2$-valued symmetric one-dimensional random walk.
We call the process $(Y_n)_{n\in\BbN}$ the \textit{vertical skeleton}.
We denote by
\[\eta_n(y)=\sum_{k=0}^n\id_{\{Y_k=y\}}, \ \ n\in\BbN, y\in\BbV_2\]
its \textit{occupation time} at level $y$.
\end{defi}

\begin{defi}{}
Suppose the vertical 
skeleton and the environments of the orientations
are given.
Let $(\xi^{(y)}_n)_{n\in\BbN^*, y\in\BbV_2}$ be a doubly infinite sequence
of independent identically distributed $\BbN$-valued 
geometric random variables of parametres $p$ and $q=1-p$. Let $(\eta_n(y))$
be the occupation times of the vertical skeleton. We call 
\textit{horizontally  embedded} random walk the process $(X_n)_{n\in\BbN}$ with
\[X_n  
=  \sum_{y\in\BbV_2} \epsilon_y \sum_{i=1}^{\eta_{n-1}(y)} \xi_i^{(y)},
\ \ n\in\BbN.\]
\end{defi}

\Rk The significance of the random variable $X_n$ is the
horizontal displacement after $n-1$ vertical moves of the skeleton $(Y_l)$.
Notice that the random walk $(X_n)$ has unbounded (although integrable)
increments. As a matter of fact, they are signed integer-valued geometric
random variables.

\begin{lemm}
Let 
\[T_n=n+\sum_{y\in\BbV_2} \sum_{i=1}^{\eta_{n-1}(y)} \xi_i^{(y)}\]
be the instant just after the random walk $(\bM_k)$ has performed
its $n^{\textrm{th}}$ vertical move (with the convention that the
sum $\sum_{i}$ vanishes whenever $\eta_{n-1}(y)=0$.) Then
\[\bM_{T_n}=(X_n,Y_n).\]
\end{lemm} 
\Proof Obvious. \eproof

Define $\sigma_0=0$ and recursively, for $n=1,2,\ldots$, 
$\sigma_n=\inf\{k\geq\sigma_{n-1}: Y_k=0\}>\sigma_{n-1}$,
the $n^{\textrm{th}}$ return to the origin for the vertical
skeleton. Then obviously, $\bM_{T_{\sigma_n}}=(X_{\sigma_n},0)$.
To study the recurrence or the transience of $(\bM_k)$, we must study
how often $\bM_k=(0,0)$. Now, $\bM_k=(0,0)$ if and only if $X_k=0$ and
$Y_k=0$. Since $(Y_k)$ is a simple random walk, the event
$\{Y_k=0\}$ is realised only at the instants $\sigma_n$, $n=0,1,2,\ldots$.

\begin{lemm}
\label{lem-return-of-M}
Let $\cF=\sigma(\psi_i, i\in\BbN)$ and $\cG=\sigma(\epsilon_y, y\in\BbV_2)$.
Denote $(\sigma_n)$ the sequence of consecutive returns to $0$ for the
skeleton random walk $(Y_k)_{k\in\BbN}$ and $Z$ a $\BbN$-valued  
random variable having the same distribution as $\xi_1$.   
Then
\[\sum_{l=0}^\infty \BbP(\bM_l=(0,0)| \cF\vee \cG)=
\sum_{n=0}^\infty \BbP(I(X_{\sigma_n},\epsilon_0 Z) \ni 0 | \cF\vee \cG),\]
where, for $x\in\BbZ$,  $z\in \BbN$, and $\epsilon=\pm1$,
$I(x,\epsilon z)=\{x, \ldots, x+z\}$ if $\epsilon =+1$ and
$\{x-z, \ldots, x\}$ if $\epsilon=-1$.
\end{lemm}
\Proof
For the process $(\bM_l)$ to return to the origin, both horizontal and
vertical components must be $0$. Since $Y_{\sigma_n}=0$ and only then     
$\bM_{T_{\sigma_n}}=(X_{\sigma_n}, 0)$. For $k=T_{\sigma_n},
\ldots, T_{\sigma_{n+1}-1}$, the
process $\bM_k$ can as a matter of fact vanish only when $k$ is in the first
part of this discrete time interval, before the process performs any vertical move,
namely if either $X_{\sigma_n}=0$ or if the points 
$X_{\sigma_n}$ and $X_{\sigma_n+1}$ stradle the point $0$.
Now
\begin{eqnarray*}
X_{\sigma_n+1}-X_{\sigma_n} &=&
\sum_{y\in\BbV_2} \epsilon_y \left(
\sum_{i=1}^{\eta_{\sigma_n}(y)} \xi_i^{(y)} -
\sum_{i=1}^{\eta_{\sigma_n-1}(y)} \xi_i^{(y)} \right)\\  
&=& \epsilon_0 \xi^{(0)}_{\eta_{\sigma_n}(0)} \\
&\elaw & \epsilon_0 Z.
\end{eqnarray*}
Therefore, the process $(\bM_l)$ can vanish for
$l\in \{T_{\sigma_n},
\ldots, T_{\sigma_{n+1}}-1 \}$ if, and only if, the point $0$ belongs    
to the set of integers $I(X_{\sigma_n},\epsilon_0 Z)$.
\eproof

\Rk Since the random variable $Z$ 
is almost surely finite (and even integrable),
the recurrence/transience properties of the random walk $(\bM_l)$
on the two-dimensional oriented lattice are essentially given by the
recurrence/transience properties of the embedded random walk
$(X_{\sigma_n})$ which is an one-dimensional random walk with unbounded
jumps in a random scenery.
Notice however that this situation is fundamentally different from the 
random walk in a random scenery studied in \cite{KesSpi}.

Although all the subsequent estimates for recurrence/transience of the process
can be carried on the right hand side of the formula obtained in lemma \ref{lem-return-of-M},
some 
can be simplified if we take advantage of the following
\begin{lemm}
\label{lem-equivalence-MX}
Let $\cF=\sigma(\psi_i, i\in\BbN)$ and $\cG=\sigma(\epsilon_y, y\in\BbV_2)$.
Denote $(\sigma_n)$ the sequence of consecutive returns to $0$ for the
skeleton random walk $(Y_k)_{k\in\BbN}$.
Then
\begin{enumerate}
\item
If $\sum_{n=0}^\infty \BbP_0(X_{\sigma_n}=0|\cF\vee \cG)=\infty$ then 
$\sum_{l=0}^\infty \BbP(\bM_l=(0,0)| \cF\vee \cG)=\infty$.
\item
If $(X_{\sigma_n})_{n\in\BbN}$ is transient then $(M_n)_{n\in\BbN}$ is also transient.
\end{enumerate}
\end{lemm}
\Proof
Notice that 
\begin{eqnarray*}
 \BbP(I(X_{\sigma_n},\epsilon_0 Z) \ni 0 | \cF\vee \cG)& = &
p\BbP(X_{\sigma_n}=0|\cF\vee \cG)\\
&&+\id_{\{\epsilon_0=-1\}}
\BbP(\cup_{x\in\BbN^*} \{ X_{\sigma_n}=x; Z\geq x|\cF\vee \cG)\\
&&+\id_{\{\epsilon_0=1\}}
\BbP(\cup_{x\in\BbN^*} \{ X_{\sigma_n}=-x; Z\geq x|\cF\vee \cG).
\end{eqnarray*}
In case 1.\ the result follows immediately from lemma \ref{lem-return-of-M}.
In case 2., since the process $(X_{\sigma_n})_{n\in\BbN}$ is transient, there exists
a constant $C>0$ such that for all $x\in\BbZ$, we have
$\sum_{n\in\BbN}\BbP_0(X_{\sigma_n}=x|\cF\vee \cG)\leq C<\infty$.
Consequently,
\[\sum_{n\in\BbN}\BbP_0(\cup_{x\in\BbN^*}\{X_{\sigma_n}=x; Z\geq x\}|\cF\vee \cG)=
\sum_{n\in\BbN}\sum_{x\geq 1} q^x \BbP_0(X_{\sigma_n}=x|\cF\vee \cG)\leq \frac{q}{1-q}C,\]
proving thus the transience of $(M_n)_{n\in\BbN}$.
\eproof

\section{Proof of theorems \ref{th-L} and \ref{th-H}}
Let $\xi$ be a geometric random variable equidistributed with $\xi_i^{(y)}$.
Denote 
\[\chi(\theta)=\BbE \exp(i\theta\xi)=\frac{p}{1-q\exp(i\theta)}
=r(\theta)\exp(i\alpha\theta), \ \ \theta
\in[-\pi,\pi]\]
its characteristic function, where
\[r(\theta)=|\chi(\theta)|=\frac{p}{\sqrt{p^2+2q(1-\cos\theta)}}=r(-\theta)\]
and
\[\alpha(\theta)=\arctan\frac{q\sin\theta}{1-q\cos\theta}=-\alpha(-\theta).\]
Notice that $r(\theta)<1$ for $\theta\in [-\pi,\\pi]\setminus\{0\}$.
Recall that we denote $\cF=\sigma(\psi_i, i\in\BbN)$ and
$\cG=\sigma(\epsilon_y, y \in\BbV_2)$. Then
\begin{eqnarray*}
\BbE\exp(i\theta X_n) &=& \BbE\left(\BbE(\exp(i\theta X_n)|\cF\vee\cG)\right)\\
&=&\BbE\left(\BbE(\exp(i\theta\sum_{y\in\BbV_2}\epsilon_y 
\sum_{i=1}^{\eta_{n-1}(y)}\xi_i^{(y)}|\cF\vee\cG)\right)\\
&=&\BbE\left(\prod_{y\in\BbV_2}\chi(\theta\epsilon_y)^{\eta_{n-1}(y)}\right).
\end{eqnarray*}

\subsection{The random walk on the $\BbL$ lattice}
\begin{lemm}
\label{lem-combinat}
For all $n\in\BbN^*$, the occupation time of the skeleton
random walk verifies:
\[\sum_{y\in\BbV_2}(-1)^y \eta_{\sigma_n-1}(y)=0.\]
\end{lemm}

\Proof
Using the strong Markov property of the process $(Y_n)$, it is enough
to show the above equality for $n=1$. Let
\[H=\sum_{y\in\BbV_2}(-1)^y \eta_{\sigma_1-1}(y).\]
Since the skeleton walk is a simple one-dimensional walk, we know
that $\sigma_1<\infty$ almost surely. Hence, 
\[H=\sum_{k\in\BbN} H\id_{\{\sigma_1=2k\}}\equiv \sum_{k\in\BbN} H_{2k}.\]
Decompose the event 
$\{\sigma_1=2k\}=\{\sigma_1=2k; Y_1=1\}\cup\{\sigma_1=2k; Y_1=-1\}$ and
consider the trajectory
\[\omega_l=\left\{\begin{array}{ll} 
l &\textrm{for}\ \ l=0,\ldots, k\\
2k-l &\textrm{for}\ \ l=k+1,\ldots, 2k.
\end{array}\right.\]
Obviously this trajectory belongs to the event
$A=\{\sigma_1=2k; Y_1=1\}$. On this trajectory,
\[
H(\omega) = 1+ \sum_{l=1}^k (-1)^l + \sum_{l=k+1}^{2k-1} (-1)^{2k-l}=0.\]
All other trajectories contributing to the event $A$ are obtained from $\omega$
by applying successively elementary transformations of the following type:
if a level $y>2$ is a local maximum of a trajectory in $A$, reflect this
local maximum with respect to the level $y-1$. The new trajectory is still
in $A$ and this operation modifies the occupation times of levels
$y$ and $y-2$ by
\begin{eqnarray*}
\eta_{\sigma_1-1}(y) &\leftarrow& \eta_{\sigma_1-1}(y)-1 \\
\eta_{\sigma_1-1}(y-2) &\leftarrow& \eta_{\sigma_1-1}(y-2)+1.
\end{eqnarray*}
The corresponding net modification in the value of $H$ is
$(-1)^y(-1)+(-1)^{y-2}(+1)=0$. The proof if completed by symmetry
for trajectories in the set 
$\{\sigma_1=2k; Y_1=-1\}$.
\eproof

\noindent
\textit{Proof of theorem \ref{th-L}:}
By lemma \ref{lem-return-of-M}, taking
expectations on both sides, we get the obvious minoration
\[\sum_{n\in\BbN}\BbP(\bM_n=(0,0))\geq \sum_{n\in\BbN}\BbP(X_{\sigma_n}=0).\]
Moreover, at every moment that the
skeleton random walk returns to the origin, the embedded random walk
starts afresh so that the process $(X_{\sigma_n})_{n\in\BbN}$ verifies
a renewal equation.
Therefore, to show recurrence of the random walk $(\bM_n)$,
it is enough to show that 
$\sum_{n\in\BbN}\BbP(X_{\sigma_n}=0)=\infty$.

Now, using parity properties of the modulus and angular part
of the characteristic function we get
\begin{eqnarray*}
\BbE\exp(i\theta X_{\sigma_n}) 
&=& \BbE\left(\prod_{y\in\BbV_2}\chi(\theta\epsilon_y)^{\eta_{\sigma_n-1}(y)}\right)\\
&=& \BbE\left(r(\theta)^{\sum_{y\in\BbV_2} \eta_{\sigma_n-1}(y)}
\exp(i\alpha(\theta)\sum_{y\in\BbV_2}\epsilon_y\eta_{\sigma_n-1}(y))\right)\\
&=& \BbE\left(r(\theta)^{\sum_{y\in\BbV_2} \eta_{\sigma_n-1}(y)}\right),
\end{eqnarray*}
by the previous combinatorial lemma \ref{lem-combinat}.
Moreover, 
\[\sum_{y\in\BbV_2} \eta_{n-1}(y)=\sum_{y\in\BbV_2}
\sum_{l=0}^{n-1} \id_{\{Y_l=y\}}=n.\]
Hence, using again the renewal equation,
$\BbE(\exp(i\theta X_{\sigma_n}))=\BbE(r(\theta)^{\sigma_n})
=(\BbE r(\theta)^{\sigma_1})^n$.
Now, $\sigma_1$ is the time of first return to the origin
for a simple random walk, its generating function reads
$f(s)=\BbE s^{\sigma_1}=
\sum_{k=1}^\infty s^{2k}\BbP(\sigma_1=2k)=1-\sqrt{1-s^2}$
for $|s|\leq1$. Hence, finally,
\[\BbE(\exp(i\theta X_{\sigma_n}))=(1-\sqrt{1-r(\theta)^2})^n,\]
so that
\[\sum_{n=0}\BbP(X_{\sigma_n}=0)=\lim\limits_{\epsilon\rightarrow 0}
2 \int_\epsilon^{\pi} \frac{1}{\sqrt{1-r(\theta)^2}} d\theta,\]
since for $\theta\in[\epsilon, \pi]$ the function $r(\theta)<1$.
Now, for $\theta\rightarrow 0$,  
$\frac{1}{\sqrt{1-r(\theta)^2}}=\cO(\frac{1}{|\theta|})$ and since
$1/|\theta|$ is a non-integrable singularity at 0,
$\sum_n \BbP(X_{\sigma_n}=0)=\infty$, proving thus the
recurrence of $(\bM_n)$.
\eproof

\subsection{The random walk on the $\BbH$ lattice}

\begin{lemm}
\label{lem-sum-epsilon-H}
Let $\epsilon_y=1$ if $y\geq 0$ and $\epsilon_y=-1$ if $y< 0$.
Denote by $(\rho_k)_{k\in\BbN}$ a sequence of independent
identically distributed Rademacher variables and $(\tau_k)_{k\in\BbN}$
a sequence of independent, identically distributed random variables,
independent of the sequence $(\rho_k)_{k\in\BbN}$, such that
$\tau_1\elaw \sigma_1$, \textit{i.e.}\ the random variables $\tau_k$
have the same law as the time of first return to the origin for the
skeleton random walk. Then
\[\sum_{y\in\BbV_2} \epsilon_y \eta_{\sigma_n-1}(y) \elaw 
\sum_{k=1}^n \rho_k (\tau_k-1) +n.\]
\end{lemm}
\Proof
We have 
\begin{eqnarray*}
\sum_{y\in\BbV_2} \epsilon_y\eta_{\sigma_n-1}(y) &=&
\sum_{y\in\BbV_2} \epsilon_y \sum_{j=1}^n \sum_{k=\sigma_{j-1}}^{\sigma_j-1}
\id_{\{Y_k=y\}}\\
&=& 
\sum_{y\ne 0} \epsilon_y \sum_{j=1}^n \sum_{k=\sigma_{j-1}+1}^{\sigma_j-1}
\id_{\{Y_k=y\}}+\epsilon_0 n.
\end{eqnarray*}
Now for every $j\in\{1,\ldots,n\}$ the process $Y_k$ has the same
sign for all 
$k\in\{\sigma_{j-1}+1, \ldots, \sigma_j-1\}$. Hence
$ \sum_{y\ne 0}\sum_{k=\sigma_{j-1}+1}^{\sigma_j-1}
\id_{\{Y_k=y\}}= \sigma_j-\sigma_{j-1}-1\elaw \tau_j-1$.
However, the contribution to the sum including the $\epsilon$ variables
 must be corrected by the sign of $Y_{\sigma_{j-1}+1}\elaw Y_1$
and since the skeleton random walk is symmetric and strongly Markovian,
we have finally 
\[\sum_{y\in\BbV_2} \epsilon_y \eta_{\sigma_n-1}(y) \elaw 
\sum_{k=1}^n \rho_k (\tau_k-1) +n.\]
\eproof

\begin{prop}
\label{prop-characteristic-embedded}
For the embedded random walk, we have
\[\BbE\left(\exp(i\theta X_{\sigma_n})\right)=
g(\theta)^n,\]
where $g(\theta)=\frac{1}{2}\chi(\theta)
\left[ \left(1-\sqrt{1-\chi(\theta)^2}\right)\exp(-i\alpha(\theta)) +
\left(1-\sqrt{1-\ol{\chi}(\theta)^2}\right)\exp(i\alpha(\theta)) \right]$.
\end{prop}
\Proof
Denote $\cD=\sigma(\rho_k, k\in\BbN)$ the $\sigma$-algebra generated
by the Rademacher variables. Then, using lemma \ref{lem-sum-epsilon-H}, 
we have,
\begin{eqnarray*}
\BbE\left(\exp(i\theta X_{\sigma_n})\right) &=&
\BbE\left(\BbE\left[\exp(i\theta X_{\sigma_n})|\cD\right]\right)\\
&=&\BbE\left(\BbE\left[r(\theta)^{\sum_{j=1}^n (\tau_j-1)}
\exp\left(i\alpha(\theta)(\sum_{j=1}^n\rho_j(\tau_j-1) +n)\right)|\cD\right]\right)\\
&=& \chi(\theta)^n \BbE\left(\BbE(r(\theta)^{\sum_{j=1}^n \tau_j}
\exp(i\alpha(\theta)\sum_{j=1}^n\rho_j(\tau_j-1) )|\cD)\right)\\
&=&  \chi(\theta)^n \BbE\left(\BbE\left[\prod_{j=1}^n( 
\chi(\theta \rho_j)^{\tau_j}
\exp(-i\alpha(\theta)\rho_j))|\cD\right]\right)\\
&=& \chi(\theta)^n\prod_{j=1}^n \BbE\left[
\left(1-\sqrt{1-\chi(\theta\rho_j)^2}\right)\exp(-i\alpha(\theta)\rho_j)\right]\\
&=& g(\theta)^n.
\end{eqnarray*}
\eproof

\begin{prop}
\label{prop-transience-H}
The random walk $(\bM_n)_{n\in\BbN}$ on the lattice $\BbH$ verifies
\[\sum_{n=0}^\infty \BbP(\bM_n=(0,0))<\infty.\]
\end{prop}
\Proof Recalling that $\epsilon_0=1$, 
from lemma \ref{lem-return-of-M} we have
\begin{eqnarray*}
\sum_{n=0}^\infty \BbP(\bM_n=(0,0)) &=& 
\sum_{n=0}^\infty \BbP(I(X_{\sigma_n},Z)\ni0)\\
&=&\sum_{n=0}^\infty \sum_{x\geq 0} \BbP(X_{\sigma_n}=-x) \BbP(Z\geq x).
\end{eqnarray*}
From the proposition \ref{prop-characteristic-embedded} we obtain
\[\BbP(X_{\sigma_n}=-x)=
\int_{-\pi}^{\pi} \exp(i\theta x) g(\theta)^n d\theta,\]
so that
\begin{eqnarray*}
\sum_{x\geq 0} \BbP(X_{\sigma_n}=-x)\BbP(Z\geq x) &=&
\int_{-\pi}^{\pi} \left(\sum_{x\geq 0} \exp(i\theta x)pq^x\right) 
g(\theta)^n d\theta\\
&=& \int_{-\pi}^{\pi} \frac{p}{1-q\exp(i\theta)} g(\theta)^nd\theta.
\end{eqnarray*}
Therefore, since $|g(\theta)|<1$ for $\theta\ne0$, 
\[\sum_{n=0}^\infty \BbP(\bM_n=(0,0))=\lim_{\epsilon\rightarrow0}
\int_\epsilon^{\pi}\left[2 \Re
\chi(\theta)\frac{1}{1-g(\theta)}\right]d\theta.\] 
Now $\lim_{\theta\rightarrow 0^+} \frac{1-g(\theta)}{\sqrt{\theta}}=\frac{1}{\sqrt{2}}$, 
therefore for $\theta\rightarrow0^+$, 
$(1-g(\theta))^{-1}=\cO(\theta^{-1/2})$ that constitutes an 
integrable singularity, while the factor
$\chi(\theta)$ does not change the singular behaviour at 0,  
proving thus the transience of the random walk.
\eproof

\section{Proof of the theorem \ref{th-O}}
The main difficulty that arises when we deal with the 
$\BbO_{\bom{\epsilon}}$
lattice stems from the fact that the embedded random walk $X_{n}$, 
which always satisfies the equation 
\[X_n=\sum_{y\in\BbV_2}\epsilon_y \sum_{i=1}^{\eta_{n-1}(y)}\xi_i^{(y)},\]
cannot any longer
 be split into independent parts when sampled on the moments
$\sigma_n$ of successive returns to the origin for the skeleton walk because
the increment $X_{\sigma_{n+1}}-X_{\sigma_n}$ is not 
independent from the increment  $X_{\sigma_{n}}-X_{\sigma_{n-1}}$,
since they may share \textit{the same} random variables $\epsilon_y$ for some
$y$. Hence the embedded random walk does not verify a renewal equation and
some new techniques are needed.

\subsection{Technical estimates}
Recall that all random variables are defined on the same probability
space $(\Omega,\cA,\BbP)$; introduce the following sub-$\sigma$-algebras:
\begin{eqnarray*}
\cH &=& \sigma(\xi_i^{(y)}, i\in\BbN, y\in\BbV_2)\\
\cG &=& \sigma(\epsilon_y,  y\in\BbV_2)\\
\cF_n &=& \sigma(\psi_i, i=1,\ldots, n),
\end{eqnarray*}
with $\cF\equiv\cF_\infty$.

Introduce the sequence of events $A_n=A_{n,1}\cap A_{n,2}$ and $B_n$ with
\begin{eqnarray*}
A_{n,1} &=& \{\omega\in\Omega: 
\max_{0\leq k\leq 2n}|Y_k|<n^{\frac{1}{2}+\delta_1}\}\ 
\ \textrm{for some}\ \ \delta_1>0,\\
A_{n,2} &=& \{\omega\in\Omega: \max_{y\in\BbV_2}\eta_{2n-1}(y)<
n^{\frac{1}{2}+\delta_2}\}\ \ \textrm{for some}\ \ \delta_2>0,\\
B_n &=& \{\omega\in A_n: \left|\sum_{y\in\BbV_2}\epsilon_y \eta_{2n-1}(y)\right|
>n^{\frac{1}{2}+\delta_3}\}\ \ \textrm{for some}\ \ \delta_3>0.
\end{eqnarray*}
Obviously $A_{n,1}, A_{n,2}$ and hence $A_n$ belong to $\cF_{2n}$;
moreover 
$B_n\subseteq A_n$ and $B_n\in \cF_{2n}\vee\cG$
. We denote in the sequel generically 
$d_{n,i}=n^{\frac{1}{2}+\delta_i}$, for $i=1,2,3$.

Since $B_n\subseteq A_n$ and both sets are $\cF_{2n}\vee\cG$-measurable,
decomposing the unity as
\[1=\id_{B_n}+\id_{A_n\setminus B_n}+ \id_{A_n^c},\]
we have
\begin{eqnarray*}
\BbP(X_{2n}=0; Y_{2n}=0|\cF\vee\cG)
&=&\id_{B_n}\id_{\{Y_{2n}=0\}}\BbP(X_{2n}=0|\cF\vee\cG)\\
&&+
\id_{A_n\setminus B_n}\id_{\{Y_{2n}=0\}}\BbP(X_{2n}= 0|\cF\vee\cG)\\
&&+
\id_{A_n^c}\id_{\{Y_{2n}=0\}}\BbP(X_{2n}= 0|\cF\vee\cG),
\end{eqnarray*}
and taking expectations on both sides of the equality, we get
\[p_n=p_{n,1}+p_{n,2}+p_{n,3},\]
where
\begin{eqnarray*}
p_n&=& \BbP(X_{2n}=0; Y_{2n}=0)\\
p_{n,1}&=& \BbP(X_{2n}=0; Y_{2n}=0;B_n)\\
p_{n,2}&=& \BbP(X_{2n}=0; Y_{2n}=0;A_n\setminus B_n)\\
p_{n,3}&=& \BbP(X_{2n}=	0; Y_{2n}=0;A_n^c).
\end{eqnarray*}

\subsection{Proof of transience of the random walk on $\BbO_{\bom{\epsilon}}$}
The transience of the random walk $(\bM_n)$ will be shown by 
establishing asymptotic estimates for the probabilities 
$p_{n,1}$, $p_{n,2}$, and $p_{n,3}$, for large $n$ showing the
summability of $p_n$ and using lemma \ref{lem-equivalence-MX} to conclude.

\begin{prop}
\label{prop-estimate-pn3}
For large $n$, there exist  $\delta>0$ and $c>0$ such that
\[p_{n,3}=\cO(\exp(-cn^\delta)).\]
\end{prop}
\Proof
Write $A_n^c=A_{n,1}^c\cup A_{n,2}^c$.

We have $\BbP(A_{n,1}^c|Y_{2n}=0)=
\BbP(\max_{0\leq k\leq 2n} |Y_k|\geq d_{n,1}|Y_{2n}=0)$
for $d_{n,1}=n^{\frac{1}{2}+\delta_1}$ and some $\delta_1>0$.
Let $a_n=[d_{n,1}]$ and $R_n=\{a_n,a_n+1,\ldots, n\}$.
With this notation,
\begin{eqnarray*}
\BbP(A_{n,1}^c|Y_{2n}=0)&=&
\sum_{y\in R_n} \BbP(\max_{0\leq k\leq 2n} |Y_k|=y|Y_{2n}=0)\\
&\leq& 2 \sum_{y\in R_n} \BbP_0(Y_{2n}=2y)\\
&=& 2 \BbP_0(Y_{2n}\geq 2a_n),
\end{eqnarray*}
by the symmetry of the skeleton random walk and the reflection
principle.
The last probability is majorised by standard methods,
\begin{eqnarray*}
\BbP(A_{n,1}^c|Y_{2n}=0)&\leq&2 \BbP_0(Y_{2n}\geq 2a_n)\\
&\leq& 2\inf_{t>0} \BbP_0(\exp(t Y_{2n})\geq \exp(ta_n))\\
&\leq& 2 \inf_{t>0} \frac{(\cosh t)^{2n}}{\exp(2ta_n)}\\
&=& 2 \exp(-\frac{a_n^2}{n})\\
&\leq& 2 \exp(-n^{2\delta_1}).
\end{eqnarray*}

In a similar way,
\[\BbP(A_{n,2}^c|Y_{2n}=0)=
\BbP(\max_{y\in\BbV_2}(\eta_{2n-1}(y)\geq d_{n,2}|Y_{2n}=0),\]
with $d_{n,2}=n^{\frac{1}{2}+\delta_2}$ for some $\delta_2>0$.
The conditional probability in the right hand side of the above
equation can be trivially majorised as
\[\BbP(\max_{y\in\BbV_2}(\eta_{2n-1}(y)\geq d_{n,2}|Y_{2n}=0)
\leq \sum_{y\in\BbV_2}\frac{\BbP(\eta_{2n-1}(y)\geq d_{n,2})}
{\BbP(Y_{2n}=0)}.\]
Now,
\[\BbP_0(\eta_{2n-1}(y)\geq d_{n,2})\leq \BbP_y(\sigma_{y,[d_{n,2}]}\leq 2n)\]
where $\sigma_{y,k}$ denotes the time of $k^{\textrm{th}}$ return to point
$y$ for the skeleton random walk. 
As a matter of fact, the occupation time, $\eta_{2n-1}(y)$,
of level $y$
 for the skeleton random walk $(Y_k)$ can exceed the threshold
$[d_{n,2}]$ whenever the random walk $(Y_k)$, starting at the origin,
attains level $y$ before time $2n-1$ and then returns to this level
 at least $[d_{n,2}]$ times before time $2n-1$. Therefore,
\begin{eqnarray*}
\BbP_0(\eta_{2n-1}(y)\geq d_{n,2})&\leq&
\BbP_y \sigma_{y,[d_{n,2}]}\leq 2n)\\
&\leq& \BbP_y(\exp(-t \sigma_{y,[d_{n,2}]})\geq \exp(-2nt))\\
&\leq& \exp(2nt)\BbE_0(\exp(-t \sigma_{[d_{n,2}]}))\\
&=&\exp(2nt)(1-\sqrt{1-\exp(-2t)})^{[d_{n,2}]}, \ \ \forall t>0\\
&\leq&\inf_{t>0} \exp(2nt)(1-\sqrt{1-\exp(-2t)})^{[d_{n,2}]}\\
&=& \exp(-cn^{\delta_2}),
\end{eqnarray*}
for some positive constant $c$, uniformly in $y$.
Using the well known estimate $\BbP_0(Y_{2n}=0)=\cO(n^{-\frac{1}{2}})$,
and the fact that the sum $\sum_{y\in\BbV_2}$ is performed on the
set $\{Y_{2n}=0\}$, containing thus at most $2n+1$ terms,
we get the overall bound
\[\BbP(A_{n,2}^c|Y_{2n}=0)\leq Cn n^{\frac{1}{2}} \exp(-cn^{\delta_2}).\]

Choosing finally $0<\delta'_2<\delta_2$ and $\delta=\min(2\delta_1,\delta'_2)>0$
we conclude that 
\[\BbP(A_{n}^c|Y_{2n}=0)\leq  2\exp(-n^{2\delta-1})+C\exp(-cn^{\delta'_2})
=\cO(\exp(-cn^\delta)).\]
\eproof

\begin{coro}
We have $\sum_{n\in\BbN}p_{n,3}<\infty$.
\end{coro}

Recall that we have 
\[X_{2n}=\sum_{y\in\BbV_2}\epsilon_y \sum_{i=1}^{\eta_{2n-1}(y)} \xi_i^{(y)}
=\sum_{k=1}^{2n} \epsilon_{Y_k}\xi_k.\]
Introduce the random variables:
\begin{eqnarray*}
N_+ &=& \sum_{k=1}^{2n} \id_{\{\epsilon_{Y_k}=1\}}\\
N_- &=& \sum_{k=1}^{2n} \id_{\{\epsilon_{Y_k}=-1\}}\\
\Delta_n &=& N_+-N_-=\sum_{y\in\BbV_2}\epsilon_y \eta_{2n-1}(y).
\end{eqnarray*}

\begin{prop}
\label{prop-estimate-pn1}
For large $n$, we have
\[p_{n,1}=\cO(\exp(-n^{\delta'}))\]
for any $\delta'\in]0,2\delta_3[$.
\end{prop}
\Proof 
Using ideas in the proof of lemma \ref{lem-equivalence-MX}, it is enough to show that
$\BbP(X_{2n}=0; Y_{2n}=0;B_n)=\cO(\exp(-n^{\delta'}))$.
Remark that
$N_+$, $N_-$, and $\Delta_n$ are $\cF_{2n}\vee\cG$-measurable and $N_++N_-=2n$.
Denoting $m_1=\BbE\xi_1$, $m_2=\BbE(\xi_1^2)$ and $s^2=m_2-m_1^2$, we have
\begin{eqnarray*}
\BbE(X_{2n}|\cF_{2n}\vee\cG)&=& m_1 \Delta_n\\
\BbE(X^2_{2n}|\cF_{2n}\vee\cG)&=& 2ns^2+m_1^2 \Delta_n^2\\
\Var(X_{2n}|\cF_{2n}\vee\cG)&=& 2ns^2.
\end{eqnarray*}
For $t\in]-\infty, -\ln q[$, define the generating function for the random variable $\xi_1$,
mamely $\phi(t)=\BbE\exp(t\xi_1)$. Obviously, for small values of $|t|$, the generating function
behaves like $\phi(t)=\exp(tm_1+t^2s^2/2+\cO(t^3))$. Hence
\[\BbE(\exp(tX_{2n})|\cF_{2n}\vee\cG)=\phi(t)^{N_+}\phi(-t)^{N_-}=
\exp(tm_1\Delta_n +t^2s^2n +\cO(t^3n)).\]
Assume for the moment that $\Delta_n>d_{n,3}$. Using Markov
inequality, we have for $t<0$,
\begin{eqnarray*}
\BbP(X_{2n}=0|\cF_{2n}\vee\cG)&\leq& \BbP(X_{2n}\leq0|\cF_{2n}\vee\cG)\\
	&\leq& \BbE(\exp(tX_{2n})|\cF_{2n}\vee\cG)\\
        &=& \exp(tm_1\Delta_n +t^2s^2n +\cO(t^3n)).
\end{eqnarray*}
Now, choose $t=-\frac{m_1n^{\delta_3-1/2}}{2 s^2}$.
Hence, on $\{\Delta_n>d_{n,3}\}$, we have
\[\BbP(X_{2n}=0|\cF_{2n}\vee\cG)\leq \exp(-\frac{m_1^2}{4s^2} n^{2\delta_3}+\cO(n^{3\delta_3-1/2})).\]

For the case $\{\Delta_n<-d_{n,3}\}$, we conclude similarly, majorising $\BbP(X_{2n}=0|\cF_{2n}\vee\cG)
\leq \exp(tm_1\Delta_n +t^2s^2n +\cO(t^3n))$, for $t\in]0,-\ln q[$, and choosing
for large $n$, $t=\frac{m_1n^{\delta_3-1/2}}{2 s^2}$.

\begin{coro}
We have
$\sum_{n\in\BbN}p_{n,1}<\infty$.
\end{coro}

To conclude about transience, it remains to estimate the probability
$p_{n,2}=\BbP(X_{2n}=0; Y_{2n}=0; A_n\setminus B_n)$.

\begin{lemm}
\label{lem-estim-conditional}
On the set $A_n\setminus B_n$, we have
\[\BbP(X_{2n}=0|\cF\vee\cG)=\cO(\sqrt{\frac{\ln n}{n}}).\]
\end{lemm} 
\Proof
Use the $\cF\vee\cG$-measurability of the variables $(\epsilon_y)_{y\in\BbV_2}$ and
$(\eta_n(y))_{y\in\BbV_2,n\in\BbN}$ to express the conditional characteristic function
of the variable $X_{2n}$ as follows:
\[\chi_1(\theta)=\BbE(\exp(i\theta X_{2n})|\cF\vee\cG)=\prod_{y\in\BbV_2}\chi(\theta\epsilon_y)^
{\eta_{2n-1}(y)}.\]
Hence,
\[\BbP(X_{2n}=0|\cF\vee\cG)=\frac{1}{2\pi}\int_{-\pi}^{\pi}
\chi_1(\theta)d\theta.\]
Now use the decomposition of $\chi$ into a  the modulus part, $r(\theta)$ --- that is an even
function of $\theta$ ---  and the angular part of $\alpha(\theta)$ 
and the fact that there is a constant $K<1$ such that
for $\theta\in[-\pi,-\pi/2]\cup[\pi/2,\pi]$ we can bound $r(\theta)<K$ to majorise
\[\BbP(X_{2n}=0|\cF\vee\cG)\leq \frac{1}{\pi}\int_{0}^{\pi/2}
r(\theta)^{2n}d\theta +\cO(K^n).\]
Fix $a_n=\sqrt{\frac{\ln n}{n}}$ and split the above integral over $[0,\pi/2]=
[0,a_n]\cup[a_n,\pi/2]$. For the first part, we majorise the integrand by 1, so that
\[\int_{0}^{a_n}      
r(\theta)^{2n}d\theta \leq a_n.\]
For the second part, use the majorisation $r(\theta)\leq \exp(-\frac{3}{8}\theta^2)$
valid for $\theta\in ]0,\pi/2]$ to estimate 
\[\frac{1}{\pi}\int_{a_n}^{\pi/2}
r(\theta)^{2n}d\theta= \cO(n^{-3/4}).\]
Since the estimate of the first part dominates, the result follows.
\eproof

\begin{prop}
\label{prop-estim-Bnc}
For all $\delta_5>0$, and for large $n$
\[\BbP(A_n\setminus B_n|\cF)=\cO(n^{-\frac{1}{4}+\delta_5}).\]
\end{prop}
\Proof
The required probability is an estimate, on the event $A_n$, of
the conditional probability $\BbP(|\sum_{y\in\BbV_2} \zeta_y|\leq d_{n,3}|\cF)$,
where we denote $\zeta_y=\epsilon_y \eta_{2n-1}(y)$. Extend the probability space
$(\Omega,\cA,\BbP)$ to carry an auxilliary variable
$G$ assumed to be centered Gaussian with variance $d_{n,3}^2$, (conditionally on $\cF$) 
independent of the $\zeta_y$'s. Since both $G$ and $\sum_{y\in\BbV_2} \zeta_y$ are
(conditionally on $\cF$) symmetric random variables and $[-d_{n,3},d_{n,3}]$
is a symmetric set around 0, then by Anderson's inequality, there exists a positive constant
$c$ such that
\[\BbP(|\sum_{y\in\BbV_2} \zeta_y|\leq d_{n,3}|\cF)\leq
c \BbP(|\sum_{y\in\BbV_2} \zeta_y+G|\leq d_{n,3}|\cF).\]
Let 
\[\chi_2(t)=\BbE(\exp(i t \sum_y \zeta_y)|\cF)=\prod_y\cos(\eta_{2n-1}(y) t),\]
and
\[\chi_3(t)=\BbE(\exp(i t G)|\cF)=\exp(-t^2d^2_{n,3}/2).\] 
Therefore, 
\[\BbE(\exp(i t (\sum_y \zeta_y+G))|\cF)=\chi_2(t)\chi_3(t),\]
and using the Plancherel's formula,
\[\BbP(|\sum_{y\in\BbV_2} \zeta_y+G|\leq d_{n,3}|\cF)=\frac{d_{n,3}}{\pi}
\int\frac{\sin(td_{n,3})}{td_{n,3}}\chi_2(t)\chi_3(t) dt\leq Cd_{n,3} I,\]
where
\[I=\int\prod_y\cos(\eta_{2n-1}(y) t)\exp(-t^2d^2_{n,3}/2) dt.\]
Fix $b_n=\frac{n^{\delta_4}}{d_{n,3}}$, for some $\delta_4>0$  and split the integral defining $I$ into
$I_1+I_2$, the first part being for $|t|\leq b_n$ and the second for
$|t|> b_n$.

We have
\begin{eqnarray*}
I_2&\leq& C \int_{|t|> b_n} \exp(-t^2 d^2_{n,3}/2) \frac{dt}{2\pi}\\
&=& \frac{C}{d_{n,3}} \int_{|s|> n^{\delta_4}} \exp(-s^2/2) \frac{ds}{2\pi}\\
&\leq&2 \frac{C}{d_{n,3}}\frac{1}{n^{\delta_4}} \frac{\exp(-n^{2\delta_4}/2)}{2\pi},
\end{eqnarray*}
because the probability that a centred normal random variable of variance 1,
whose density is denoted $\phi$,  exceeds a threshold
$x>0$ is majorised by $\frac{\phi(x)}{x}$.

For $I_1$ we get,
\[ I_1\leq \int_{|t|\leq b_n} \prod_y |\cos(\eta_{2n-1}(y)t)| dt.\]
Now, $\sum_y\frac{\eta_{2n-1}(y)}{2n}=1$. Therefore, applying Hölder's
inequality we obtain
\begin{eqnarray*}
I_1 &\leq& \prod_y \left[
\left( 
\int_{|t|\leq b_n} |\cos(\eta_{2n-1}(y)t)|^{\frac{2n}{\eta_{2n-1}(y)}} dt
\right)^{\frac{\eta_{2n-1}(y)}{2n}}\right]\\
&\leq& \sup_{y: \eta_{2n-1}(y)\ne 0}\int_{|t|\leq b_n} 
|\cos(\eta_{2n-1}(y)t)|^{\frac{2n}{\eta_{2n-1}(y)}} dt,
\end{eqnarray*}
because the terms in $|cos(\cdot)|$ in the integrand are less than 1 and for $x\in[0,1]$ and 
$p\geq1$ we have that $x^p\leq x$. 

Now, on the set $A_n$ and  for every $y:\eta_{2n-1}(y)\ne 0$,   we have
$|\eta_{2n-1}(y)\ne 0 |\leq b_n d_{n,2}=\frac{d_{n,2}}{d_{n,3}}n^{\delta_4}$
and we can always choose the parametres $\delta_2, \delta_3, \delta_4$ so that
$\delta_2+\delta_4-\delta_3=-\delta_5<0$. For those $y$,
\begin{eqnarray*}
I_1 &\leq& 2b_n \int_{|v|\leq \pi/2} 
|\cos(\frac{\eta_{2n-1}(y)}{d_{n,3}}n^{\delta_4} v)|^{\frac{2n}{\eta_{2n-1}(y)}} dv\\
&\leq& 2b_n \int_{|v|\leq \pi/2} |\cos v|^{\frac{2n}{\eta_{2n-1}(y)}} dv.
\end{eqnarray*}
For $|v|< \pi/2$ we have that $|\cos v|\leq \exp(-cv^2/2)$, so that
\begin{eqnarray*}
I_1&\leq& 2 b_n \int \exp(-c\frac{2nv^2}{\eta_{2n-1}(y)}) dv\\
&\leq& cb_n \frac{\eta_{2n-1}(y)^{1/2}}{n^{1/2}} \int \exp(-v^2/2) \frac{dv}{\sqrt 2\pi}\\
&\leq& c\frac{d_{n,2}^{1/2}}{d_{n,3}} \frac{n^{\delta_4}}{n^{1/2}},
\end{eqnarray*}
and finally the overall probability is majorised by 
$\cO(d_{n,3}(I_1+I_2))=\cO(n^{-1/4})$.
\eproof

\begin{coro}
\[\sum_{n\in\BbN}p_{n,2}<\infty.\]
\end{coro}
\Proof
Recalling that for the standard random walk 
$\BbP(Y_{2n}=0)=\cO(n^{-1/2})$ and from the estimates obtained in
\ref{lem-estim-conditional} and \ref{prop-estim-Bnc}, we have 
\begin{eqnarray*}
p_{n,2}&=&\BbP(X_{2n}=0; Y_{2n}=0;A_n\setminus B_n)\\
&=& \BbE(\BbE\left(\id_{Y_{2n}=0}\left[
\BbE(\id_{A_n\setminus B_n} \BbP(X_{2n}=0|\cF\vee\cG) | \cF)\right]\right))\\
&=& \cO( n^{-1/2} n^{-1/4} \sqrt{\frac{\ln n}{n}})\\
&=& \cO(n^{-5/4} \ln n),
\end{eqnarray*}
proving thus the summability of $p_{n,2}$.
\eproof

\begin{theo}
For almost all realisations of the random environment $\bom{\epsilon}$, 
the random walk on the lattice $\BbO_{\bom{\epsilon}}$
is transient.
\end{theo}
\Proof
The transience  is a simple consequence
the previous propositions. As a matter of fact
$p_n=p_{n,1}+p_{n,2}+p_{n,3}$ is summable because the partial probabilities
$p_{n,i}$, for $i=1,2,3$ are all shown to be summable.
\eproof

\subsection{Strong law of large numbers for the random walk on $\BbO_{\bom{\epsilon}}$}

\begin{prop}
\label{prop-w-LLN}
Denote $m_1=\BbE\xi_1$. Then conditionally on $\cF\vee\cG$, the ratio
$\frac{X_n-m_1\Delta_n}{n}$ converges in probability to 0, \textit{i.e.}\
for all $\lambda>0$, 
\[\lim_{n\rightarrow\infty} \BbP(|\frac{X_n-m_1\Delta_n}{n}|\geq \lambda|\cF\vee\cG)=0.\]
\end{prop}
\Proof
Since $X_n=\sum_y \epsilon_y \sum_{i=1}^{\eta_{n-1}(y)} \xi_i^{(y)}$, it follows that
$\BbE(X_n|\cF\vee\cG)= m_1(N_+-N_-)=m_1\Delta_n$.
Letting $m_2=\BbE(\xi^2_1)$ and $s^2=m_2-m_1^2$, we have by developing
\begin{eqnarray*}
X_n^2 &=& \sum_y \sum_{i=1}^{\eta_{n-1}(y)} (\xi_i^{(y)})^2 +
\sum_y \sum_{i_1=1}^{\eta_{n-1}(y)}\sum_{i_2=1}^{\eta_{n-1}(y)} \xi_{i_1}^{(y)}\xi_{i_2}^{(y)}\\
&& + \sum_{y_1} \sum_{y_2\ne y_1} \epsilon_{y_1}\epsilon_{y_2}
\sum_{i_1=1}^{\eta_{n-1}(y_1)}\sum_{i_2=1}^{\eta_{n-1}(y_2)} \xi_{i_1}^{(y_1)}\xi_{i_2}^{(y_2)},
\end{eqnarray*}
that $\BbE(X^2_n|\cF\vee\cG)=s^2 n + m_1^2 \Delta_n^2$. Consequently,
$\BbE((X_n-m_1\Delta_n)^2|\cF\vee\cG)=s^2 n$ and the result follows by a
straightforward application of Chebychev's inequality. 
\eproof

\begin{prop}
We have 
\[\lim_{n\rightarrow\infty} \frac{\Delta_n}{n}= 0\ \ \textrm{almost surely}.\]
\end{prop}
\Proof
Since $\Delta_n=\sum_y\epsilon_y \eta_{n-1}(y)$, using the symmetry of the random 
variables $\epsilon_y$, it follows immediately that
$\BbE(\Delta_n|\cF)=0$.
For some positive integer $r$, compute $\BbE(\Delta_n^{2r}|\cF)$. Using again the symmetry
of the random variables $\epsilon_y$, only terms containing even powers of
$\epsilon_y$ will survive. Among these terms, the dominant one for large $n$ is the term
$\sum_{y_1} \eta_{n-1}(y_1)^2 \ldots \sum_{y_r} \eta_{n-1}(y_r)^2$ and
each sum appearing is estimated by
\begin{eqnarray*}
\sum_{y} \eta_{n-1}(y)^2&\leq& \max_y \eta_{n-1}(y) \sum_{y} \eta_{n-1}(y)\\
&=& \max_y \eta_{n-1}(y) n\\
&\leq& d_{n,2} n \id_{A_{n,2}}+ n^2 \id_{A_{n,2}^c}.
\end{eqnarray*}
It follows that
\[\BbE(\Delta_{n}^{2r}|\cF)\leq n^r d_{n,3}^r \id_{A_{n,2}}+ n^{2r}\id_{A_{n,2}^c}\]
and consequently, using the estimate $\BbP(A^c_{n,2})\leq \exp(-n^{\delta_2})$, obtained
in the proof of proposition \ref{prop-estimate-pn3}, we get
\[\BbE(\Delta_n^{2r})\leq n^{3r/2+r\delta_2}+n^{2r}\exp(-n^{\delta_2}),\]
so that, for $r\geq 3$, 
\[\sum_{n\in\BbN} \BbE(\frac{\Delta_n^{2r}}{n^{2r}})< \infty.\]
The result follows by straightforward application of Borel-Cantelli lemma.
\eproof

\Rk Using similar arguments (choosing $r$ sufficiently large), 
it can easily be shown that for every $\beta>3/4$,
\[\lim_{n\rightarrow \infty} \frac{\Delta_n}{n^\beta}=0\ \ \textrm{almost surely.}\]

\begin{theo}
The embedded random walk $(X_n)$ on the $\BbO_{\bom{\epsilon}}$ lattice has almost surely
 zero speed, \textit{i.e.}\
\[\lim_{n\rightarrow\infty} \frac{X_n}{n}=0, \ \ \textrm{almost surely.}\]
\end{theo}
\Proof
It is enough to show that for some positive integer $r$, 
\[\sum_{n\in\BbN} \BbE(\frac{(X_n-m_1\Delta_n)^{2r}}{n^{2r}})< \infty\]
since then the almost sure convergence to 0 of $X_n/n$ will follow from
the almost sure convergence of $\Delta_n/n$ to 0.
Following exactly the same scheme as in the previous proposition,
start by developing  $\BbE(X_n^{2r}|\cF\vee\cH)$. The symmetry of the random
variables $\epsilon_y$ over which we integrate, guarantees that only
terms with even powers of each random variable $\epsilon_y$ will remain.
Perform now the integration over the random variables $\xi$.  As it was
the case in the proof of the previous proposition, the expectation 
$\BbE(X_n^{2r}|\cF)$ is estimated
by $C(\sum_y\eta_{n-1}(y)^2)^r$ and the individual sum over $y$ is again
estimated as previously. Thus choosing $s\geq 3$ the result follows from
Borel-Cantelli lemma. 
\eproof
  
\section{Conclusion, open problems, and further developments}
It is shown that random walks on oriented lattices exhibit novel and interesting
features and arise in many situations in topological field theories. It was quite surprising
for us to discover that so little was previously known about this kind of random walk, in sharp
contrast with random walks on undirected lattices for which there are several hundreds of
papers and also excellent books (like \cite{Spi,Woe} not to mention but the two more
complete ones.)

For the alternate lattice $\BbL$, several different proofs
for the recurrence can be given; we chose to present here the most 
elementary one. The main point 
is that the lattice is still periodic in the vertical direction. 
So group arguments can be used to prove the recurrence instead of our proof.
Another possibility is to regroup even and odd ordinates to
a new effective lattice that it is not any more oriented. This approach has been
used in \cite{Guillotin} for Manhattan lattices. 
It is easy to see that instead of taking period 2 in the vertical direction we use another
arbitrary periodicity, \textit{i.e.}\ take horizontal strips of width $l$ such that all
the horizontal lines inside a given strip are going in the same direction but the
directions alternate for every strip to the following one, the lattice is fundamentally
the same as $\BbL$ so that the random walk is still recurrent.
The overall asymptotic behaviour of all these random walks is quite reminiscent of
the random walk on unoriented ordinary $\BbZ^2$ lattice.

The lattice $\BbH$ is quite different. Here the random walk is transient. 
It can be shown that it has still a zero speed but a non zero angular
speed, so that asymptotically it has infinite winding number around the origin.
The study of the detailed behaviour of this walk is postponed in a subsequent publication.

The lattice $\BbO_{\bom{\epsilon}}$ has a very rich structure. 
The reason for which this lattice is transient is the presence and the size of
fluctuations. It should be interesting to ask whether other random environments
modify this characteristic.
Another interesting question is whether the  random walk on this
lattice verifies a functional limit theorem possibly with unconventional normalisation.
All these questions are under investigation.

\bibliographystyle{plain}
\scriptsize{
\bibliography{rwre}
 }            

\end{document}